\documentclass[11pt]{amsart}
\usepackage{amsmath}
\usepackage{amssymb}
\usepackage{amsfonts}
\usepackage{amsthm}
\usepackage{enumerate}
\usepackage[mathscr]{eucal}
\usepackage{verbatim}
\usepackage{amsthm}
\usepackage{amscd}
\usepackage[mathscr]{eucal}
\usepackage{appendix}
\usepackage{tikz}
\usepackage{stmaryrd}

\numberwithin{equation}{section}

\newtheorem{innercustomgeneric}{\customgenericname}
\providecommand{\customgenericname}{}

\newcommand{\newcustomtheorem}[2]{\newenvironment{#1}[1]
  {\renewcommand\customgenericname{#2}
   \renewcommand\theinnercustomgeneric{##1}\innercustomgeneric}{\endinnercustomgeneric}}

\newcustomtheorem{customthm}{Theorem}

\newcommand{\newcustomlemma}[2]{\newenvironment{#1}[1]
  {\renewcommand\customgenericname{#2}
   \renewcommand\theinnercustomgeneric{##1} \innercustomgeneric}{\endinnercustomgeneric}}

\newcustomlemma{customlemma}{Lemma}

\newcustomlemma{customproposition}{Proposition}

\newcustomlemma{customclaim}{Claim}

\theoremstyle{plain}
\newtheorem{theorem}{Theorem}[section]
\newtheorem{lemma}[theorem]{Lemma}
\newtheorem{corollary}[theorem]{Corollary}
\newtheorem{proposition}[theorem]{Proposition}

\theoremstyle{remark}

\theoremstyle{definition}
\newtheorem{defn}{Definition}

\newcommand{\q}{\quad}
\newcommand{\qq}{\qquad}

\newcommand{\bbz}{\mathbb{Z}}

\newcommand{\bbr}{\mathbb{R}}

\newcommand{\bbrn}{\mathbb R^n}

\newcommand{\bbn}{\mathbb{N}}

\newcommand{\QQ}{\mathcal{Q}}

\newcommand{\xxxi}{\vec{\boldsymbol{\xi}\,}}
\newcommand{\xxx}{\vec{\boldsymbol{x}}}

\newcommand{\ppp}{\vec{\boldsymbol{p}}}

\newcommand{\uuu}{\vec{\boldsymbol{u}}}

\def\000{\vec{\boldsymbol{0}}}

\newcommand{\supp}{\mathrm{supp}}

\newcommand{\wh}{\widehat}

\newcounter{question}

\newcommand{\bpf}{\begin{proof}}

\newcommand{\epf}{\end{proof}}

\allowdisplaybreaks

\makeatletter
\@namedef{subjclassname@2020}{\textup{2020} Mathematics Subject Classification}
\makeatother

\makeindex         

\begin{document}

\author{Bae Jun Park}
\address{B. Park, Department of Mathematics, Sungkyunkwan University, Suwon 16419, Republic of Korea}
\email{bpark43@skku.edu}

\author{Naohito Tomita}
\address{N. Tomita, Department of Mathematics, Graduate School of Science, Osaka University, Toyonaka, Osaka 560-0043, Japan}
\email{tomita@math.sci.osaka-u.ac.jp}

\thanks{B. Park is supported in part by NRF grant 2022R1F1A1063637 and by POSCO Science Fellowship of POSCO TJ Park Foundation. B. Park is grateful for support by the Open KIAS Center at Korea Institute
for Advanced Study. N. Tomita was supported by JSPS KAKENHI Grant Number 20K03700.}

\title[Multilinear pseudo-differential operators of type $(0,0)$]{Sharp Maximal function estimates for multilinear pseudo-differential operators of type $(0,0)$}

\begin{abstract} 

In this paper, we study sharp maximal function estimates for multilinear pseudo-differential operators.
Our target is operators of type $(0,0)$ for which a differentiation does not make any decay of the associated symbol. 
Analogous results for operators of type $(\rho,\rho)$, $0<\rho<1$, appeared in an earlier work of the authors \cite{Park_Tomita_submitted}, but a different approach is given for $\rho=0$.
\end{abstract}

\maketitle


\section{Introduction}
 
  Let $l$ be a positive integer greater or equal to $1$ in this paper.
 For $0<r<\infty$ and  locally $r$-th integrable functions $f_1,\dots,f_l$ on $\bbrn$, we  define the  $l$-sublinear  Hardy-Littlewood maximal function by
$$\mathbf{M}_r(f_1,\dots,f_l)(x):=\sup_{Q:x\in Q}\bigg( \frac{1}{|Q|^l}\int_{Q\times\cdots\times Q}\prod_{j=1}^{l}\big|f_j(u_j) \big|^r\; d\uuu \bigg)^{\frac{1}{r}}$$
where $d\uuu:=du_1\cdots du_l$ and the supremum is taken over all cubes containing $x$.
For a locally $r$th integrable function $f$ on $\bbrn$, we also define the (inhomogeneous and sublinear) sharp maximal function by
 \begin{align*}
 \mathcal{M}_r^{\sharp}f(x)&:=\sup_{Q:x\in Q, \ell(Q)\ge 1}\Big(\frac{1}{|Q|}\int_Q\big| f(y)\big|^r \; dy\Big)^{\frac{1}{r}}\\
 &\q\q+\sup_{Q:x\in Q, \ell(Q)<1}\inf_{c_Q\in \mathbb{C}}\Big(\frac{1}{|Q|}\int_Q \big| f(y)-c_Q\big|^r\; dy\Big)^{\frac{1}{r}}
 \end{align*}
where the supremum ranges over all cubes in $\bbrn$ containing $x$ whose side-length is $\ge 1$ in the first one, and is $< 1$ in the second one.
Then 
\begin{equation}\label{bmoproperty}
\Vert f\Vert_{bmo(\bbrn)}:=\big\Vert \mathcal{M}_1^{\sharp}f     \big\Vert_{L^{\infty}(\bbrn)}\sim \big\Vert \mathcal{M}_r^{\sharp}f     \big\Vert_{L^{\infty}(\bbrn)}
\end{equation} for $1<r<\infty$,
where $bmo$ denotes the localized version of the usual space of functions with bounded mean oscillation.

\hfill

 Given $m\in\bbr$ and $0\le \rho<1$, the $l$-linear H\"ormander symbol class $\mathbb{M}_lS_{\rho,\rho}^m(\bbrn)$ consists of infinitely many differentiable functions $\sigma(x,\xi_1,\dots,\xi_l)$ on $(\bbrn)^{l+1}$ satisfying the property that for any multi-indices $\alpha, \beta_1,\dots,\beta_l\in \{ 0,1,2,\dots\}^n$ there exists a constant $C_{\alpha,\beta_1,\dots,\beta_l}>0$ such that
$$\big| \partial_{x}^{\alpha}\partial_{\xi_1}^{\beta_1}\cdots\partial_{\beta_l}^{\beta_l}\sigma(x,\xi_1,\dots,\beta_l)\big|\le C_{\alpha,\beta_1,\dots,\beta_l}\big( 1+|\xi_1|+\cdots+|\beta_l|\big)^{m-\rho(|\beta_1|+\dots+|\beta_l|-|\alpha|)}.$$
Let $\mathscr{S}(\bbrn)$ denote the Schwartz class on $\bbrn$.
 For a symbol $\sigma\in \mathbb{M}_lS_{\rho,\rho}^{m}(\bbrn)$, the corresponding $l$-linear pseudo-differential operator $T_{\sigma}$ is defined as
 $$T_{\sigma}\big(f_1,\dots,f_l\big)(x):=\int_{(\bbrn)^l}\sigma(x,\xi_1,\dots,\xi_l)\prod_{j=1}^{l}\wh{f_j}(\xi_j) \, e^{2\pi i\langle x,\xi_1+\dots+\xi_l\rangle}\; d\xxxi$$
 for $f_1,\dots,f_l\in \mathscr{S}(\bbrn)$, where $\wh{f}(\xi):=\int_{\bbrn} f(x)e^{-2\pi i\langle x,\xi\rangle}\;dx$ is the Fourier transform of $f$, and $d\xxxi:=d\xi_1\cdots d\xi_l$.
 Denote by $\mathrm{Op}\mathbb{M}_lS_{\rho,\rho}^{m}(\bbrn)$ the family of such operators with symbols in $\mathbb{M}_lS_{\rho,\rho}^{m}(\bbrn)$.

Then we have the following  sharp maximal function estimates which originated in Chanillo and Torchinsky \cite{Ch_To1985} for linear operators. See also Naibo \cite{Na2015} for bilinear ones.
\begin{customthm}{A}\cite{Park_Tomita_submitted}\label{thma}
Let $1<r\le 2$, $0<\rho<1$, and $l\ge 1$. Suppose that $m\le -\frac{nl}{r}(1-\rho)$ and $\sigma\in \mathbb{M}_lS_{\rho,\rho}^{m}(\bbrn)$. 
\begin{enumerate}
\item If $l=1$, then
$$\mathcal{M}_1^{\sharp}\big( T_{\sigma}f\big)(x)\lesssim \mathrm{M}_rf(x), \q \q x\in\bbrn$$
for $f\in \mathscr{S}(\bbrn)$.
\item If $l\ge 2$, then
 \begin{equation*}
 \mathcal{M}^{\sharp}_{\frac{r}{l}}\big( T_{\sigma}(f_1,\dots,f_l)\big)(x)\lesssim \mathbf{M}_r\big(f_1,\dots,f_l\big)(x), \q \q ~x\in\bbrn
 \end{equation*}
for all $f_1,\dots,f_l\in\mathscr{S}(\bbrn)$.
\end{enumerate}
\end{customthm}
 Here and in the sequel, the symbol $A\lesssim B$ indicates that $A\le CB$ for some constant $C>0$ independent of the variable quantities $A$ and $B$, and $A\sim B$ if $A\lesssim B$ and $B\lesssim A$ hold simultaneously.

The key idea to prove Theorem \ref{thma} is the kernel estimates in \cite[Lemmas 2.1 and 3.1]{Park_Tomita_submitted}.
In order to illustrate this, take a simple case when $\sigma(x,\xi)=\sigma(\xi)$ in $\mathbb{M}_1S_{\rho,\rho}^m(\bbrn)$.
 By applying (inhomogeneous) Littlewood-Paley decomposition, which will be actually discussed in the next section,
we can write $T_{\sigma}f=\sum_{k\in\bbn_0}{T_{\sigma_k}f}$ and the kernel $K_k(x,y)$ of the operator $T_{\sigma_k}$ satisfies the size estimate
\begin{equation}\label{keyideaprevious}
|K_k(x,y)|\lesssim_M 2^{-k(\rho N-m-n)}\frac{1}{|x-y|^N}, \qquad |x-y|\gtrsim 1
\end{equation} by using integration by parts $N$ times. If $\rho>0$, by taking $N$ sufficiently large, we have enough exponential decay $2^{-k(\rho N-m-n)}$ (even though $\rho$ is very small).  Essentially similar properties are inherent in \cite[Lemmas 2.1 and 3.1]{Park_Tomita_submitted}. 
This is not that serious when $k$ is small enough, but since we will finally take the infinite sum over $k\in \{0,1,2,\dots\}$, appropriate decay in $k$ should be necessary to make it summable. However, when $\rho=0$ the decay property is no longer available and this prevents Theorem \ref{thma} from being extended to $\rho=0$.

\hfill

In this paper, we will establish analogous pointwise estimates for $\rho=0$, overcoming the obstacle mentioned above regarding \eqref{keyideaprevious}.
Our main result is
\begin{theorem}\label{mainthm}
Let $1<r\le 2$ and $l\ge 1$. Suppose that $m\le -\frac{nl}{r}$ and $\sigma\in \mathbb{M}_lS_{0,0}^{m}(\bbrn)$. 
Then we have
 \begin{equation*}
 \mathcal{M}^{\sharp}_{\frac{r}{l}}\big( T_{\sigma}(f_1,\dots,f_l)\big)(x)\lesssim \mathbf{M}_r\big(f_1,\dots,f_l\big)(x), \q \q ~x\in\bbrn
 \end{equation*}
for all $f_1,\dots,f_l\in\mathscr{S}(\bbrn)$.
\end{theorem}
As an extension of Theorem \ref{thma} to $\rho=0$, Theorem \ref{mainthm}  will be proved in a seemingly similar structure. However, we would like to emphasize that if $\rho=0$, 
sufficient exponential decay for $k$ cannot be derived through multiple integration by parts, as seen in \eqref{keyideaprevious}, so a different approach that can compensate for the critical issue is needed. For this one, we will adopt a nontrivial decomposition of a function $f$ (e.g.  \eqref{decompfep2} and  \eqref{dildecfg} below).\\

Now we recall multi-variable extension of $A_p$ weight classes, introduced by Lerner, Ombrosi, P\'erez, Torres, and Trujillo-Gonz\'alez \cite{Le_Om_Pe_To_Tr2009}.
 \begin{defn}
 Let $1< p_1,\dots,p_l<\infty$ and $1/p=1/p_1+\dots+1/p_l$.
 Then we define $\mathrm{A}_{\ppp}$, $\ppp=(p_1,\dots,p_l)$, to be the class of $l$-tuples of weights $\vec{w}:=(w_1,\dots,w_l)$ satisfying
 $$\sup_Q\bigg[ \Big( \frac{1}{|Q|}\int_Q v_{\vec{w}}(x)\; dx \Big)^{1/p}\prod_{j=1}^{l}\Big(\frac{1}{|Q|}\int_Q \big(w_j(x)\big)^{1-p_j'}\; dx \Big)^{1/p_j}\bigg]<\infty$$ 
 where $p_j'$ denotes the H\"older conjugate of $p_j$ and
 $$v_{\vec{w}}(x):=\prod_{j=1}^{l}\big( w_j(x) \big)^{p/p_j}.$$
 \end{defn}
If $l=1$, the class $A_{\vec{p}}$ coincides with the classical Muckenhoupt $A_p$ class.
It is known in \cite{Le_Om_Pe_To_Tr2009} that for $1<p_1,\dots,p_l<\infty$ with $1/p=1/p_1+\cdots+1/p_l$,
\begin{equation}\label{weightcha}
\vec{w}=(w_1,\dots,w_l)\in A_{\ppp}~\text{ if and only if }~\big\Vert \mathbf{M}(f_1,\dots,f_l)\big\Vert_{L^p(v_{\vec{w}})}\lesssim \prod_{j=1}^{l}\Vert f_j\Vert_{L^{p_j}(w_j)}
\end{equation}
 for all locally integrable functions $f_1,\dots,f_l$,  where $\mathbf{M}:=\mathbf{M}_1$.
As discussed in \cite[Theorem 1.6]{Park_Tomita_submitted}, Theorem \ref{thma}, together with \eqref{weightcha}, deduces a weighted norm inequality for $T_\sigma\in \mathrm{Op}\mathbb{M}_lS_{\rho,\rho}^{-\frac{nl}{r}(1-\rho)}(\bbrn)$, involving $\vec{w}=(w_1,\dots, w_l)\in A_{\frac{p_1}{r},\dots,\frac{p_l}{r}}$, and it can be naturally extended to $\rho=0$ replacing Theorem \ref{thma} by Theorem \ref{mainthm}.
  \begin{theorem}\label{weightedthm}
 Let $l\ge 1$, $1< r\le 2$, and $r<p_1,\dots,p_l<\infty$ with $1/p=1/p_1+\dots+1/p_l$. 
 Suppose that  $\sigma\in \mathbb{M}_lS_{0,0}^{-\frac{nl}{r}}(\bbrn)$.
If an $l$-tuple of weights $\vec{w}=(w_1,\dots, w_l)$ belongs to the class $A_{\frac{p_1}{r},\dots,\frac{p_l}{r}}$, then the weighted norm inequality
$$\big\Vert T_{\sigma}(f_1,\dots,f_l)\big\Vert_{L^p(v_{\vec{w}})}\lesssim \prod_{j=1}^{l}\Vert f_j\Vert_{L^{p_j}(w_j)}$$
holds 
for $f_1,\dots,f_l\in \mathscr{S}(\bbrn)$.
 \end{theorem}
 The proof is exactly same as that of \cite[Theorem 1.6]{Park_Tomita_submitted}, simply employing Theorem \ref{mainthm} instead of Theorem \ref{thma}.
 So it is omitted here. 

We also have the following end-point estimate as a corollary of Theorem \ref{mainthm}.
\begin{corollary}
Suppose that
$\sigma\in \mathbb{M}_lS_{0,0}^{-\frac{nl}{2}}(\bbrn)$.
Then we have
\begin{equation}\label{infbmoest}
\big\Vert T_{\sigma}(f_1,\dots,f_l)\big\Vert_{BMO(\bbrn)}\lesssim \prod_{j=1}^{l}\Vert f_j\Vert_{L^{\infty}(\bbrn)}
\end{equation}
for all $f_1,\dots,f_l\in\mathscr{S}(\bbrn)$.
\end{corollary}
The inequality \eqref{infbmoest} follows from the $BMO$ characterization
$$\Vert f\Vert_{BMO(\bbrn)}\sim_p \sup_{Q}\inf_{c_Q}
\left( \frac{1}{|Q|}\int_Q |f(x)-c_Q|^p\, dx \right)^{1/p}, \qq 0<p<\infty$$
and Theorem \ref{mainthm}. See \cite[Corollary 1.5]{Park_Tomita_submitted} for more details.

\hfill

{\bf Notation.}
Let $\bbn$ and $\bbz$ be the sets of all natural numbers and integers, respectively, and $\bbn_0:=\bbn\cup \{0\}$.
For each set $U$ in $\bbrn$  we denote by $\chi_U$ its characteristic function and by $U^c$ its complement.
For each cube $Q$ in $\bbrn$, let $\ell(Q)$ denote the side-length of $Q$.

\section{Preliminaries}

\subsection{Decomposition of Pseudo-differential operators}
 
 Let $\Phi$ be a Schwartz function on $\bbr^{nl}$ such that its Fourier transform $\wh{\Phi}$ is equal to $1$ on the unit ball centered at the origin and is supported in the ball of radius $2$.
 Let $\Psi\in\mathscr{S}(\bbr^{nl})$ satisfy $\wh{\Psi}(\xxxi):=\wh{\Phi}(\xxxi)-\wh{\Phi}(2\xxxi)$ for $\xxxi\in \bbr^{nl}$.
 For each $k\in\bbz$, we define $\Psi_k(\xxx):=2^{knl}\Psi(2^k \xxx)$ for $\xxx\in \bbr^{nl}$. Then $\{\Phi\}\cup \{\Psi_k\}_{k\in\bbn}$ forms inhomogeneous Littlewood-Paley partition of unity. Note that
 $$\supp(\wh{\Psi_k})\subset \big\{\xxxi\in \bbr^{nl} : 2^{k-1}\le |\xxxi|\le 2^{k+1}\big\}, \q k\in\bbz$$
 and
 \begin{equation}\label{linearLPdecom}
 \wh{\Phi}(\xxxi)+\sum_{k\in\bbn}\wh{\Psi_k}(\xxxi)=1.
\end{equation}
Using \eqref{linearLPdecom}, we decompose  $\sigma \in \mathbb{M}_lS_{0,0}^{m}(\bbrn)$ as
$$\sigma(x,\xxxi)=\sigma(x,\xxxi)\wh{\Phi}(\xxxi)+\sum_{k\in\bbn}\sigma(x,\xxxi)\wh{\Psi_k}(\xxxi)=:\sigma_0(x,\xxxi)+\sum_{k\in\bbn}\sigma_k(x,\xxxi), \q \xxxi\in\bbr^{nl}$$
 and then express
\begin{equation*}
 T_{\sigma}(f_1,\dots,f_l)=\sum_{k\in\bbn_0}T_{\sigma_k}(f_1,\dots,f_l)
\end{equation*}
 where $T_{\sigma_k}$ are the $l$-linear pseudo-differential operators associated with $\sigma_k\in \mathbb{M}_lS_{0,0}^{m}(\bbrn)$. Here, we notice that for all multi-indices $\alpha,\beta_1,\dots,\beta_l \in (\bbn_0)^n$
 \begin{equation*}
 \big| \partial_x^{\alpha}\partial_{\xi_1}^{\beta_1}\cdots\partial_{\xi_l}^{\beta_l}\sigma_k(x,\xi_1,\dots,\xi_l)\big|\lesssim_{\alpha,\beta_1,\dots,\beta_l}\big( 1+|\xi_1|+\cdots+|\xi_l|\big)^{m}
 \end{equation*}
 where the implicit constant in the inequality is independent of $k$. In such a case, we will say 
 $\sigma_k$ belongs to $\mathbb{M}_lS_{0,0}^{m}(\bbrn)$ uniformly in $k$.

\hfill

\subsection{Generalized trace theorem for Sobolev spaces}

For $s>0$ and $N\in\bbn$, let $(I_N-\Delta_N)^{s/2}$ denote the Bessel potential acting on  tempered distributions on $\bbr^N$, defined by
$$(I_N-\Delta_N)^{\frac{s}{2}}F:= \big( (1+4\pi^2  |\cdot|^2)^{\frac{s}{2}}\wh{F} \,\big)^{\vee}.$$
The Sobolev space $L^2_s(\mathbb{R}^N)$ consists of all $F \in L^2(\mathbb{R}^N)$ satisfying
$$\Vert F \Vert_{L^2_s(\bbr^N)}:=\big\Vert (I_N-\Delta_N)^{\frac{s}{2}}F \big\Vert_{L^2(\bbr^N)}=\Big(\int_{\bbr^N}\big(1+ 4\pi^2|\xi|^2 \big)^s\big|\wh{F}(\xi) \big|^2\; d\xi \Big)^{\frac{1}{2}}<\infty.$$

 \begin{lemma}\cite[Lemma 5.1]{Park_Tomita_submitted}\label{tracethm}
Let $l\ge 2$ and $s>0$.
 Then
 \begin{equation*}
 \Vert \widetilde{G} \Vert_{L^2_{s}(\bbrn)}
\lesssim_{s,r,n} \Vert G\Vert_{L^2_{s+\frac{(l-1)n}{2}}((\bbrn)^l)}
 \end{equation*}
 for Schwartz functions $G$ on $(\bbrn)^l$,
 where $\widetilde{G}(x)=G(x_1,\dots,x_{l}) |_{x_1=\dots=x_l=x}$ and $x,x_1,\dots,x_l \in \bbrn$.
 \end{lemma}

\hfill

\subsection{General $l$-linear estimates for $T_{\sigma}\in \mathbb{M}_lS_{0,0}^{m}(\bbrn)$}
  Let $h^p$ denote the local Hardy space, introduced by Goldberg \cite{Go1979}. The boundedness of (linear) pseudo-differential operators of type $(0,0)$ on $h^p$ was given by Calder\'on and Vaillancourt \cite{Ca_Va1972} for $p=2$, Coifman and Meyer \cite{Co_Me1978} for $1<p<\infty$, and Miyachi \cite{Mi1987} and P\"aiv\"arinta and Somersalo \cite{Pa_So1988} for $0<p \le 1$. The multilinear setting was studied by B\'enyi and Torres \cite{Be_To2004}, Michalowski, Rule, and Staubach \cite{Mi_Ru_St2014}, B\'enyi, Bernicot, Maldonado, Naibo, and Torres \cite{Be_Be_Ma_Na_To2013}, Miyachi and Tomita \cite{Mi_To2013}, and finally Kato, Miyachi, and Tomita \cite{Ka_Mi_To2022} proved the following.
  \begin{customthm}{B}\label{thme}
 Let $l\ge 1$ and $0<p_1,\dots,p_l\le \infty$ with $1/p_1+\dots+1/p_l=1/p$. Suppose that $m\in\bbr$ and $\sigma\in \mathbb{M}_lS_{0,0}^m(\bbrn)$.
 Then 
  $$m\le m_0(\ppp):=-n\bigg[ \sum_{j=1}^{l}\max\Big(\frac{1}{p_j},\frac{1}{2} \Big) -\min\Big(\frac{1}{p},\frac{1}{2} \Big)\bigg], \q \ppp:=(p_1,\dots,p_l)$$
 if and only if
 $$\big\Vert T_{\sigma}\big(f_1,\dots,f_l\big)\big\Vert_{Y^p(\bbrn)}\lesssim \prod_{j=1}^{l}\Vert f_j\Vert_{Y^{p_j}(\bbrn)}$$
 for $f_1,\dots,f_l\in \mathscr{S}(\bbrn)$,
 where 
 $$Y^p:=\begin{cases}
 h^p & 0<p<\infty\\
 bmo & p=\infty
 \end{cases}.$$
  \end{customthm}
Note that the inequality $\|f\|_{L^{p}} \lesssim \|f\|_{h^{p}}$ holds for locally integrable functions $f$, where $0<p \le 1$.
Thus Theorem \ref{thme} deduces that if $1<r \le 2$ and $\sigma\in \mathbb{M}_lS_{0,0}^{-n(\frac{l}{r}-\frac{1}{2})}(\bbrn)$, then
\begin{equation}\label{lrtolrlbd}
\big\Vert T_{\sigma}\big(f_1,\dots,f_l\big)\big\Vert_{L^{r/l}(\bbrn)}\lesssim \prod_{j=1}^{l}\Vert f_j\Vert_{L^{r}(\bbrn)} .
\end{equation}

\section{Key Estimates}

We fix $k\in \bbn_0$.
 Let $P_k$ be the concentric dilate of $Q$ with $\ell(P_k)\ge 10\sqrt{n} l \ell(Q)$ and $\chi_{P_k}$ denote its characteristic function.
We write
 \begin{align*}
 f_{j}= f_{j}\chi_{P_k}+f_j\chi_{P_k^c}=:f_{j,k}^{(0)}+f_{j,k}^{(1)} \q \text{ for }~j=1,\dots,l
 \end{align*}
 so that each $T_{\sigma_k}(f_1,\dots,f_l)$ can be decomposed as 
 $$T_{\sigma_k}(f_1,\dots,f_l)=\sum_{\mu_1,\dots,\mu_l\in \{0,1\}}T_{\sigma_k}\big(f_{1,k}^{(\mu_1)},f_{2,k}^{(\mu_2)},\dots,f_{l,k}^{(\mu_l)}\big).$$

Then a key estimate for the proof of Theorem \ref{mainthm} is 
 \begin{proposition}\label{keypropo1}
Let $1< r\le 2$, $l\ge 1$, $k\in\bbn_0$, and $\ell(Q)<1$. 
 Suppose that  $\sigma\in \mathbb{M}_lS_{0,0}^{-\frac{nl}{r}}(\bbrn)$ and $x\in Q$.
Then we have
 \begin{equation}\label{keypropo1ineq}
 \inf_{c_Q\in \mathbb{C}}\bigg(\frac{1}{|Q|}\int_Q \big| T_{\sigma_k}\big(f_{1,k}^{(0)},\dots,f^{(0)}_{l,k} \big)(y)-c_Q\big|^{\frac{r}{l}} \; dy\bigg)^{\frac{l}{r}}\lesssim \ell(P_k)^{\frac{nl}{r}}\mathbf{M}_r\big(f_1,\dots,f_l\big)(x)
 \end{equation}
uniformly in $k$.
 \end{proposition}

 In order to prove Proposition \ref{keypropo1}, we first need the following lemma, which is actually an extension of \cite[Lemma 5.4]{Park_Tomita_submitted} to $\rho=0$.
    \begin{lemma}\label{Naibopropo}
Let $1<r\le 2$ and $l\ge 1$. Suppose that $\sigma\in \mathbb{M}_lS_{0,0}^{-\frac{nl}{r}}(\bbrn)$. Then we have
 $$\big\Vert T_{\sigma}(f_1,\dots,f_l) \big\Vert_{bmo(\bbrn)}\lesssim \prod_{j=1}^{l}\Vert f_j\Vert_{L^r(\bbrn)}$$
 for $f_1,\dots,f_l\in \mathscr{S}(\bbrn)$.
 \end{lemma} 
 \begin{proof}
  The linear case $l=1$ has been already proved by Miyachi and Yabuta \cite[Lemma 5.6]{Mi_Ya1987} by using duality argument with $(h^1)^*=bmo$, but we provide another proof, using Sobolev embeddings into $bmo$.
 We notice that the operator
 $$  (I_n-\Delta_n)^{\frac{n}{4}}T_{\sigma} (I_n-\Delta_n)^{\frac{1}{2}n(\frac{1}{r}-\frac{1}{2})}$$
 is a pseudo-differential operator associated with a symbol in $\mathbb{M}_1S_{0,0}^{0}(\bbr^{n})$.
Then its $L^2(\bbr^{n})$ boundedness, in view of Theorem \ref{thme}, deduces that
 \begin{align*}
 \big\Vert T_{\sigma}f\big\Vert_{L^2_{\frac{n}{2}}(\bbr^{n})}
 &=\Big\Vert \Big(  (I_n-\Delta_n)^{\frac{n}{4}}T_{\sigma} (I_n-\Delta_n)^{\frac{1}{2}n(\frac{1}{r}-\frac{1}{2})}    \Big)  (I_n-\Delta_n)^{-\frac{1}{2}n(\frac{1}{r}-\frac{1}{2})}f   \Big\Vert_{L^2(\bbr^{n})}\\
 &\lesssim   \Big\Vert  (I_n-\Delta_n)^{-\frac{1}{2}n(\frac{1}{r}-\frac{1}{2})}f\Big\Vert_{L^2(\bbr^{n})}              =\big\Vert f\big\Vert_{L^2_{-n(\frac{1}{r}-\frac{1}{2})}(\bbr^{n})}\lesssim \Vert f\Vert_{L^r(\bbr^{n})}
 \end{align*}
 for any Schwartz functions $f$ on $\bbr^{n}$,
 where the last inequality follows from the Sobolev embedding theorem, but we note that $r>1$.
 Therefore, by the embedding $L_{\frac{n}{2}}^2(\bbrn)\hookrightarrow bmo(\bbrn)$ (e.g. \cite[Theorem 2.3]{Park_2019Nachr}),
 we obtain
 \begin{align*}
 \big\Vert T_{\sigma}f\big\Vert_{bmo(\bbrn)}\lesssim  \big\Vert T_{\sigma}f\big\Vert_{L^2_{\frac{n}{2}}(\bbr^{n})}\lesssim \Vert f\Vert_{L^r(\bbrn)},
 \end{align*}
as desired.

 Now we assume $l\ge 2$. In this case, we mimic the proof of \cite[Lemma 5.4]{Park_Tomita_submitted}.
 We first define the linear symbol $\Sigma$ in $\bbr^{nl}\times \bbr^{nl}$ as
 \begin{equation}\label{Sigmadef}
 \Sigma(X,\Xi):=\sigma(x_1,\xi_1,\dots,\xi_l),\q X:=(x_1,\dots,x_l)\in (\bbr^{n})^l, ~ \Xi:=(\xi_1,\dots,\xi_l)\in (\bbr^{n})^l.
 \end{equation}
 Then it is obvious that the symbol $\Sigma$ belongs to the linear H\"ormander class $\mathbb{M}_1S_{0,0}^{-\frac{nl}{r}}(\bbr^{nl})$.
Moreover, the operator
 $$  (I_{nl}-\Delta_{nl})^{\frac{nl}{4}}T_{\Sigma} (I_{nl}-\Delta_{nl})^{\frac{1}{2}nl(\frac{1}{r}-\frac{1}{2})}$$
 is a linear pseudo-differential operator in $\mathrm{Op}\mathbb{M}_1S_{0,0}^{0}(\bbr^{nl})$.
Now, by applying its $L^2(\bbr^{nl})$ boundedness in  Theorem \ref{thme}, we have
 \begin{align*}
 \big\Vert T_{\Sigma}(G)\big\Vert_{L^2_{\frac{nl}{2}}(\bbr^{nl})}
 &=\Big\Vert \Big( (I_{nl}-\Delta_{nl})^{\frac{nl}{4}}T_{\Sigma}(I_{nl}-\Delta_{nl})^{\frac{1}{2}nl(\frac{1}{r}-\frac{1}{2})}    \Big)  (I_{nl}-\Delta_{nl})^{-\frac{1}{2}nl(\frac{1}{r}-\frac{1}{2})}G   \Big\Vert_{L^2(\bbr^{nl})}\\
 &\lesssim   \Big\Vert  (I_{nl}-\Delta_{nl})^{-\frac{1}{2}nl(\frac{1}{r}-\frac{1}{2})}G\Big\Vert_{L^2(\bbr^{nl})}              =\big\Vert G\big\Vert_{L^2_{-nl(\frac{1}{r}-\frac{1}{2})}(\bbr^{nl})}\lesssim \Vert G\Vert_{L^r(\bbr^{nl})}
 \end{align*}
 for any Schwartz functions $G$ on $\bbr^{nl}$,
 where the Sobolev embedding theorem is applied in the last inequality.
Finally,   since $$ T_{\sigma}\big(f_1,\dots, f_l \big)(x)=T_{\Sigma}(f_1\otimes\cdots\otimes f_l)(x_1,\dots,x_l) |_{x_1=\dots=x_l=x}, \q x\in\bbrn ,$$
the embedding $L_{\frac{n}{2}}^2(\bbrn)\hookrightarrow bmo(\bbrn)$ and Lemma \ref{tracethm} with $s=\frac{n}{2}$ deduce
 \begin{align*}
 \big\Vert T_{\sigma}(f_1,\dots,f_l)\big\Vert_{bmo(\bbrn)}&\lesssim \big\Vert T_{\sigma}(f_1,\dots,f_l)\big\Vert_{L_{\frac{n}{2}}^2(\bbrn)}\\
 &\lesssim  \big\Vert T_{\Sigma}(f_1\otimes\cdots\otimes f_l)\big\Vert_{L^2_{\frac{nl}{2}}(\bbr^{nl})}\\
 &\lesssim \Vert f_1\otimes\cdots\otimes f_l\Vert_{L^r(\bbr^{nl})}=\prod_{j=1}^{l}\Vert f_j\Vert_{L^r(\bbrn)}.
 \end{align*}
This completes the proof.
 \end{proof}

 Now we prove Propositions \ref{keypropo1}. 
 \begin{proof}[Proof of Proposition \ref{keypropo1}]
In view of \eqref{bmoproperty}, together with H\"older's inequality if $l\ge 2$,
the left-hand side of \eqref{keypropo1ineq} is no more than
$$\Big\Vert T_{\sigma_k}\big( f_{1,k}^{(0)},\dots,f_{l,k}^{(0)}\big)\Big\Vert_{bmo(\bbrn)}$$
and this is, via Lemma \ref{Naibopropo}, bounded by a constant times
 \begin{align}\label{l2onpest}
\prod_{j=1}^{l}\Vert f_j\Vert_{L^r(P_k)}&=\Big( \int_{P_k\times\cdots\times P_k} \big| f_1(u_1)\cdots f_l(u_l)\big|^{r}\; du_1\cdots du_l\Big)^{\frac{1}{r}}\nonumber\\
&\lesssim \ell(P_k)^{\frac{ nl}{r}}\mathbf{M}_r\big(f_1,\dots,f_l\big)(x)
 \end{align} 
 as $x\in Q\subset P_k$.
 Then \eqref{keypropo1ineq} follows.
 \end{proof}

If one of $\mu_j$ is equal to $1$, we benefit from the following proposition.
 \begin{proposition}\cite[Propositions 2.3 and 5.3]{Park_Tomita_submitted}\label{keypropo2}
 Let $1\le r\le 2$, $l\ge 1$, $k\in\bbn_0$, and $x,y\in Q$. Suppose that $\sigma\in \mathbb{M}_lS_{0,0}^{-\frac{nl}{r}}(\bbrn)$. If each $\mu_j$ is $0$ or $1$, excluding $\mu_1=\cdots=\mu_l=0$, then we have
 \begin{equation}\label{keypropo2_1}
 \big| T_{\sigma_k}\big(f_{1,k}^{(\mu_1)},\dots,f_{l,k}^{(\mu_l)}\big)(y)\big|\lesssim_N \ell(P_k)^{-(N-\frac{nl}{r})}\mathbf{M}_r\big(f_1,\dots,f_l\big)(x)
\end{equation}
and
\begin{align}\label{keypropo2_2}
& \big| T_{\sigma_k}\big(f_{1,k}^{(\mu_1)},\dots,f_{l,k}^{(\mu_l)}\big)(y)-T_{\sigma_k}\big(f_{1,k}^{(\mu_1)},\dots,f_{l,k}^{(\mu_l)}\big)(x)\big|\nonumber\\
&\lesssim_N 2^k\ell(Q)\ell(P_k)^{-(N-\frac{nl}{r})}\mathbf{M}_r\big(f_1,\dots,f_l\big)(x)
\end{align}
for any $N>\frac{nl}{r}$.
\end{proposition}

\hfill

Another key estimate is
\begin{proposition}\label{newkeypropo}
 Let $1<r\le 2$, $l\ge 1$, $x\in Q$, and $\ell(Q)<1$.
 Suppose that $\sigma\in \mathbb{M}_lS_{0,0}^{-\frac{nl}{r}}(\bbrn)$.
 Then we have
 \begin{equation}\label{1qaveragemm}
 \bigg( \frac{1}{|Q|}\int_Q  \Big| \sum_{k:2^k\ell(Q)\ge 1} T_{\sigma_k}\big(f_1,\dots,f_l \big)(y)\Big|^{\frac{r}{l}}\; dy\bigg)^{\frac{l}{r}}\lesssim  \mathbf{M}_r\big(f_1,\dots,f_l \big)(x).
 \end{equation}
 \end{proposition}

\begin{proof}
We note that $2^k\ell(Q)\ge 1$ implies $k\ge 1$. Therefore, we don't need to consider the case $k=0$.
 We choose
 $
 0<\epsilon< \frac{r}{2l}( \le 1 )
$
and let $\QQ_{k,\epsilon}^{*}$ be the concentric dilate of $Q$ with $\ell(\QQ_{k,\epsilon}^{*})=10\sqrt{n}l \big( 2^k\ell(Q)\big)^{\epsilon} (\ge 10\sqrt{n}l>10\sqrt{n}l\ell(Q))$.
 Then we write
 \begin{equation}\label{decompfep2}
 f_j=f_j\chi_{\QQ_{k,\epsilon}^{*}}+f_j\chi_{(\QQ_{k,\epsilon}^{*})^c}=:\mathrm{f}_{j,k}^{(0)}+\mathrm{f}^{(1)}_{j,k} \q \text{ for }~ j=1,\dots,l
 \end{equation}
 so that the left-hand side of \eqref{1qaveragemm} is less than the sum of
 \begin{equation*}
 \mathcal{J}_1:=\bigg(  \sum_{k:2^k\ell(Q)\ge 1} \bigg( \frac{1}{|Q|}\int_Q \Big| T_{\sigma_k}\big(\mathrm{f}^{(0)}_{1,k},\dots,\mathrm{f}^{(0)}_{l,k} \big)(y)      \Big|^{\frac{r}{l}}\; dy \bigg)^{\min{(1,\frac{l}{r})}} \bigg)^{\max{(1,\frac{l}{r})}} 
 \end{equation*}
 and
 \begin{equation*}
 \mathcal{J}_2:=\sum_{\substack{\mu_1,\dots,\mu_l\in \{0,1\}\\ (\mu_1,\dots,\mu_l)\not= (0,\dots,0) }}\bigg(\sum_{k:2^k\ell(Q)\ge 1} \bigg( \frac{1}{|Q|}\int_Q   \Big| T_{\sigma_k}\big(\mathrm{f}^{(\mu_1)}_{k,1},\dots,\mathrm{f}^{(\mu_l)}_{k,l}\big)(y)      \Big|^{\frac{r}{l}}\; dy \bigg)^{\min{(1,\frac{l}{r})}}   \bigg)^{\max{(1,\frac{l}{r})}}
 \end{equation*}
where the subadditivities of $\Vert \cdot\Vert_{\ell^{r}}$ and $\Vert \cdot \Vert_{\ell^{\frac{r}{l}}}^{\frac{r}{l}}$ if $l\ge 2$ are applied.

We consider $\mathcal{J}_2$ first.
 By applying \eqref{keypropo2_1} with $N>\frac{nl}{r}$, if $(\mu_1,\dots,\mu_l)\in \{0,1\}^l\setminus \{(0,\dots,0)\}$,
 \begin{align*}
 &\bigg( \frac{1}{|Q|}\int_Q \Big| T_{\sigma_k}\big(\mathrm{f}^{(\mu_1)}_{k,1},\dots,\mathrm{f}^{(\mu_l)}_{k,l}\big)(y)      \Big|^{\frac{r}{l}}\; dy\bigg)^{\frac{l}{r}}\\
 & \lesssim \ell(\QQ_{k,\epsilon}^{*})^{-(N-\frac{nl}{r})} \mathbf{M}_r\big(f_1,\dots,f_l \big)(x)\sim \big( 2^k\ell(Q)\big)^{-\epsilon(N-\frac{nl}{r})} \mathbf{M}_r\big(f_1,\dots,f_l \big)(x)
 \end{align*}
 so that
 \begin{align*}
 \mathcal{J}_2&\lesssim \mathbf{M}_r\big(f_1,\dots,f_l \big)(x) \bigg(\sum_{k:2^k\ell(Q)\ge 1}   \big(2^k \ell(Q) \big)^{-\epsilon(\frac{rN}{l}-n)\min{(1,\frac{l}{r})}}     \bigg)^{\max{(1,\frac{l}{r})}}\\
 &\lesssim \mathbf{M}_r\big(f_1,\dots,f_l \big)(x).
 \end{align*}

 Now it remains to estimate $\mathcal{J}_1$
 and we will separately treat the cases $l=1$ and $l\ge 2$.
 
Let us first consider the case $l=1$.
We define
 $\gamma_{k}:= \Psi_{k-1}+\Psi_k+\Psi_{k+1}\in \mathscr{S}(\bbrn)$  so that
 $$\Psi_k\ast \gamma_k=\Psi_k ~ \text{ for }~k\ge 1.$$
This yields 
 $$T_{\sigma_k}\mathrm{f}^{(0)}_k=T_{\sigma_k}\big(\gamma_k \ast \mathrm{f}^{(0)}_k \big).$$
 By applying H\"older's inequality  and observing $2^{\frac{kn}{r}}\sigma_k \in \mathbb{M}_1S_{0,0}^{0}(\bbrn)$ uniformly in $k$,
 the $L^2$ boundedness for $2^{\frac{kn}{r}}T_{\sigma_k}$ in Theorem \ref{thme} deduces
 \begin{align*}
\bigg( \frac{1}{|Q|}\int_Q \big| T_{\sigma_k}\mathrm{f}^{(0)}_k(y)      \big|^r\; dy\bigg)^{\frac{1}{r}}&\le \ell(Q)^{-\frac{n}{2}}\Big\Vert T_{\sigma_k}\big(\gamma_k \ast \mathrm{f}^{(0)}_k \big)\Big\Vert_{L^2(\bbrn)}\\
&\lesssim \ell(Q)^{-\frac{n}{2}}2^{-\frac{kn}{r}}\big\Vert \gamma_k \ast \mathrm{f}^{(0)}_k  \big\Vert_{L^2(\bbrn)}
 \end{align*}
and this is bounded by
\begin{align*}
\ell(Q)^{-\frac{n}{2}}2^{-\frac{kn}{r}}\Vert \gamma_k\Vert_{L^q(\bbrn)}\big\Vert \mathrm{f}^{(0)}_k \big\Vert_{L^r(\bbrn)}\sim \big(2^k\ell(Q)\big)^{-\frac{n}{2}}\big\Vert f \big\Vert_{L^r(\QQ_{k,\epsilon}^*)}\lesssim  \big(2^k\ell(Q)\big)^{-n(\frac{1}{2}-\frac{\epsilon}{r})}\mathrm{M}_rf(x)
\end{align*}
where we applied Young's convolution inequality with $1+\frac{1}{2}=\frac{1}{q}+\frac{1}{r}$.
This finally proves
$$ \mathcal{J}_1\lesssim   \mathrm{M}_rf(x) \sum_{k:2^k\ell(Q)\ge 1} \big(2^k\ell(Q)\big)^{-n(\frac{1}{2}-\frac{\epsilon}{r})}      \lesssim \mathrm{M}_rf(x),$$
as desired.

Now assume $l\ge 2$. In this case, we will employ the trace theorem in Lemma \ref{tracethm} and the $L^2(\bbr^{nl})$ boundedness of a linear pseudo-differential operator belonging to $\mathbb{M}_1S_{0,0}^{0}(\bbr^{nl})$.
Since $2<\frac{r}{\epsilon l}$, we may choose $t$ such that
 \begin{equation*}
 2<t<\frac{r}{\epsilon l}.
 \end{equation*}
 Then by using H\"older's inequality with $\frac{r}{l}<t$ and the Sobolev embedding theorem,
 \begin{align*}
 \bigg(  \frac{1}{|Q|}\int_Q \Big| T_{\sigma_k}\big(\mathrm{f}^{(0)}_{1,k},\dots,\mathrm{f}^{(0)}_{l,k} \big)(y)      \Big|^{\frac{r}{l}}\; dy \bigg)^{\frac{l}{r}}&\le \frac{1}{|Q|^{\frac{1}{t}}}    \Big\Vert T_{\sigma_k}\big(\mathrm{f}^{(0)}_{1,k},\dots,\mathrm{f}^{(0)}_{l,k} \big)\Big\Vert_{L^t(\bbrn)}\\
 &\lesssim  \ell(Q)^{-\frac{n}{t}}\Big\Vert T_{\sigma_k}\big(\mathrm{f}^{(0)}_{1,k},\dots,\mathrm{f}^{(0)}_{l,k} \big)\Big\Vert_{L_{\frac{n}{2}-\frac{n}{t}}^2(\bbrn)}.
 \end{align*}
 Let  $\Gamma_{k}:= \Psi_{k-1}+\Psi_k+\Psi_{k+1}\in \mathscr{S}(\bbr^{nl})$ so that $\Psi_k\ast \Gamma_k=\Psi_k$ for $k\ge 1$.
Setting $\Sigma$ as in \eqref{Sigmadef} and $\Sigma_k:=\Sigma \cdot \wh{\Psi_k}$, we write
$$ T_{\sigma_k}\big(\mathrm{f}^{(0)}_{1,k},\dots,\mathrm{f}^{(0)}_{l,k} \big)(y)=T_{\Sigma_k}F_k(y,\dots,y), \q y\in\bbrn$$
 where $T_{\Sigma_k}\in \mathrm{Op}\mathbb{M}_1S_{0,0}^{-\frac{nl}{r}}(\bbr^{nl})$ and 
 $F_k:=\Gamma_k \ast\big( \mathrm{f}^{(0)}_{1,k}\otimes \cdots \otimes\mathrm{f}^{(0)}_{l,k}\big)$ is a function defined on $\bbr^{nl}$.
 Therefore, applying Lemma \ref{tracethm} with $\frac{n}{2}-\frac{n}{t}>0$, 
 \begin{align*}
 \Big\Vert T_{\sigma_k}\big(\mathrm{f}^{(0)}_{1,k},\dots,\mathrm{f}^{(0)}_{l,k} \big)\Big\Vert_{L_{\frac{n}{2}-\frac{n}{t}}^2(\bbrn)}&\lesssim \Big\Vert T_{\Sigma_k}F_k\Big\Vert_{L_{\frac{nl}{2}-\frac{n}{t}}^2(\bbr^{nl})}\\
 &=\Big\Vert \Big(     (I_{nl}-\Delta_{nl})^{\frac{1}{2}(\frac{nl}{2}-\frac{n}{t})}T_{\Sigma_k}     \Big)       F_k\Big\Vert_{L^2(\bbr^{nl})}.
 \end{align*}
Observing
$$ 2^{kn(\frac{l}{r}+\frac{1}{t}-\frac{l}{2})} (I_{nl}-\Delta_{nl})^{\frac{1}{2}(\frac{nl}{2}-\frac{n}{t})}T_{\Sigma_k} \in \mathrm{Op}\mathbb{M}_1S_{0,0}^{0}(\bbr^{nl}) \q \text{ uniformly in }~k,  $$
 its $L^2(\bbr^{nl})$ boundedness yields that the $L^2$ norm in the last expression is bounded by a constant times
 \begin{align*}
 2^{-kn(\frac{l}{r}+\frac{1}{t}-\frac{l}{2})}\big\Vert   F_k\big\Vert_{L^2(\bbr^{nl})}&\lesssim 2^{-kn(\frac{l}{r}+\frac{1}{t}-\frac{l}{2})}\big\Vert \Gamma_k\Vert_{L^q(\bbr^{nl})}\big\Vert     \mathrm{f}^{(0)}_{1,k}\otimes \cdots \otimes\mathrm{f}^{(0)}_{l,k}   \big\Vert_{L^r(\bbr^{nl})}\\
 &\sim 2^{-kn(\frac{l}{r}+\frac{1}{t}-\frac{l}{2})}2^{knl(1-\frac{1}{q})}      \prod_{j=1}^{l} \big\Vert  f_j  \big\Vert_{L^r(\QQ_{k,\epsilon}^{*})}=2^{-\frac{kn}{t}}  \prod_{j=1}^{l} \big\Vert  f_j  \big\Vert_{L^r(\QQ_{k,\epsilon}^{*})}
 \end{align*}
 where the inequality follows from Young's inequality with  $1+\frac{1}{2}=\frac{1}{q}+\frac{1}{r}$.
 Using \eqref{l2onpest}, the last displayed expression is controlled by a constant multiple of 
 $$2^{-\frac{kn}{t}}    \big( 2^k\ell(Q)\big)^{\frac{\epsilon nl}{r}}\mathbf{M}_r\big(f_1,\dots,f_l \big)(x).$$
 Finally, we have proved that
 \begin{align*}
 & \bigg(  \frac{1}{|Q|}\int_Q \Big| T_{\sigma_k}\big(\mathrm{f}^{(0)}_{1,k},\dots,\mathrm{f}^{(0)}_{l,k} \big)(y)      \Big|^{\frac{r}{l}}\; dy \bigg)^{\frac{l}{r}}\lesssim \big(2^k\ell(Q) \big)^{-n(\frac{1}{t}-\frac{\epsilon l}{r})}\mathbf{M}_r\big(f_1,\dots,f_l \big)(x)
 \end{align*}
 and this clearly implies that
 \begin{align*}
\mathcal{J}_1 &\lesssim \mathbf{M}_r\big(f_1,\dots,f_l \big)(x)\bigg( \sum_{k:2^k\ell(Q)\ge 1}       \big(2^k\ell(Q) \big)^{-n(\frac{r}{tl}-\epsilon)}        \bigg)^{\frac{l}{r}}\lesssim \mathbf{M}_r\big(f_1,\dots,f_l \big)(x).
 \end{align*}
 This completes the proof.
   \end{proof}


\section{Proof of Theorem \ref{mainthm}}

Fixing a cube $Q$ in $\bbrn$ and a point $x\in Q$, our focus shifts to demonstrating
\begin{equation}\label{mainestimate0}
\Big(\frac{1}{|Q|}\int_Q  \big|T_{\sigma}\big(f_1,\dots,f_l\big)(y)\big|^{\frac{r}{l}} \; dy \Big)^{\frac{l}{r}}\lesssim \mathbf{M}_r\big(f_1,\dots,f_l \big)(x) \q \text{ if }~\ell(Q)\ge 1
\end{equation}
and
\begin{equation}\label{mainestimate}
\inf_{c_Q\in\mathbb{C}}\Big(\frac{1}{|Q|}\int_Q  \big|T_{\sigma}\big(f_1,\dots,f_l\big)(y)-c_Q\big|^{\frac{r}{l}} \; dy \Big)^{\frac{l}{r}}\lesssim \mathbf{M}_r\big(f_1,\dots,f_l \big)(x) \q \text{ if }~ \ell(Q)<1
\end{equation}
 uniformly in $Q$ and $x$.

\subsection{Proof of \eqref{mainestimate0}} 
 Assume $\ell(Q)\ge 1$.
 Let $0<\epsilon<\frac{r}{2l}$. 
 For each $k\in\bbn_0$, 
 let $Q_{k,\epsilon}^{*}$ denote the concentric dilate of $Q$ with $\ell(Q_{k,\epsilon}^{*})=10\sqrt{n} l2^{\epsilon k}\ell(Q)$.
Then we write
  \begin{equation}\label{dildecfg}
 f_j=f_j\chi_{Q^{*}_{k,\epsilon}}+f_j\chi_{(Q^{*}_{k,\epsilon})^c}=:\mathbf{f}_{j,k}^{(0)}+\mathbf{f}_{j,k}^{(1)} \q \text{ for }~j=1,\dots,l
 \end{equation}
 so that $T_{\sigma}(f_1,\dots,f_l)$ can be decomposed as
 \begin{align*}
 T_{\sigma}(f_1,\dots,f_l)=\sum_{k\in\bbn_0}\sum_{\mu_1,\dots,\mu_l\in \{0,1\}}T_{\sigma_k}\big(\mathbf{f}_{1,k}^{(\mu_1)},\dots,\mathbf{f}_{l,k}^{(\mu_l)} \big).
 \end{align*}
Therefore, the left-hand side of \eqref{mainestimate0} is controlled by a constant times
$$\sum_{\mu_1,\dots,\mu_l\in \{0,1\}}\bigg( \sum_{k\in\bbn_0}\bigg( \frac{1}{|Q|}\int_Q \big|   T_{\sigma_k}\big(\mathbf{f}_{1,k}^{(\mu_1)},\dots,\mathbf{f}_{l,k}^{(\mu_l)}\big)(y)     \big|^{\frac{r}{l}}\; dy\bigg)^{\min{(1,\frac{l}{r})}}\bigg)^{\max(1,\frac{l}{r})}.$$

Now we claim that for each $k\in\bbn_0$,
\begin{equation}\label{llinearclaim1}
 \bigg(\frac{1}{|Q|}\int_Q \big|T_{\sigma_k}\big(\mathbf{f}_{1,k}^{(0)},\dots,\mathbf{f}_{l,k}^{(0)}\big)(y) \big|^{\frac{r}{l}} \; dy \bigg)^{\frac{l}{r}}\lesssim 2^{-kn (\frac{1}{2}-\frac{\epsilon l}{r})} \mathbf{M}_r\big(f_1,\dots,f_l\big)(x)
\end{equation}
and for $N>\frac{nl}{r}$
\begin{equation}\label{llinearclaim2}
 \bigg(\frac{1}{|Q|}\int_Q \big|T_{\sigma_k}\big(\mathbf{f}_{1,k}^{(\mu_1)},\dots,\mathbf{f}_{l,k}^{(\mu_l)}\big)(y) \big|^{\frac{r}{l}} \; dy \bigg)^{\frac{l}{r}}\lesssim 2^{-\epsilon k(N-\frac{nl}{r})}\mathbf{M}_r\big(f_1,\dots,f_l\big)(x)
\end{equation}
if at least one of $\mu_j$ is equal to $1$.
 Then \eqref{mainestimate0}  follows obviously.

 By observing
 $$2^{\frac{kn}{2}}\sigma_k\in \mathbb{M}_lS_{0,0}^{-n(\frac{l}{r}-\frac{1}{2})}(\bbrn) \q \text{ uniformly in }~k$$
 and applying \eqref{lrtolrlbd},
 the left-hand side of \eqref{llinearclaim1} is bounded by
 \begin{align*}
\frac{1}{|Q|^{\frac{l}{r}}}\Big\Vert T_{\sigma_k}\big(\mathbf{f}^{(0)}_{1,k},\dots,\mathbf{f}_{l,k}^{(0)}\big)\Big\Vert_{L^{\frac{r}{l}}(\bbrn)} \lesssim 2^{-\frac{kn}{2}}\ell(Q)^{-\frac{nl}{r}}\prod_{j=1}^{l}\Vert f_j\Vert_{L^r(Q_{k,\epsilon}^{*})}.
 \end{align*}
Using the argument in \eqref{l2onpest} with $P_k=Q_{k,\epsilon}^{*}$, we have
$$\prod_{j=1}^{l}\Vert f_j\Vert_{L^r(Q_{k,\epsilon}^{*})}\lesssim \ell(Q_{k,\epsilon}^{*})^{\frac{nl}{r}}\mathbf{M}_r\big(f_1,\dots,f_l\big)(x)\sim \big(2^{\epsilon k}\ell(Q)\big)^{\frac{nl}{r}} \mathbf{M}_r\big(f_1,\dots,f_l\big)(x)$$
 and this yields \eqref{llinearclaim1}.
 
 Now assume that one of $\mu_j$ is $1$.
 By employing \eqref{keypropo2_1},
 for  $N>\frac{nl}{r}$ and  $y\in Q$, we have
  \begin{align*}
 \Big| T_{\sigma_k}\big(\mathbf{f}_{1,k}^{(\mu_1)},\dots,\mathbf{f}_{l,k}^{(\mu_l)} \big)(y)\Big|&\lesssim  \big(\ell(Q_{k,\epsilon}^{*})\big)^{-(N-\frac{nl}{r})}\mathbf{M}_r\big(f_1,\dots,f_l\big)(x)\\
 &\sim \big(2^{\epsilon k }\ell(Q)\big)^{-(N-\frac{nl}{r})}\mathbf{M}_r\big(f_1,\dots,f_l\big)(x)\\
 &\le 2^{-\epsilon k(N-\frac{nl}{r})}\mathbf{M}_r\big(f_1,\dots,f_l\big)(x)
 \end{align*} 
 where the last inequality holds due to $\ell(Q)\ge 1$.
 This completes the proof of  \eqref{llinearclaim2}.

 \subsection{Proof of \eqref{mainestimate}}
 Assume $\ell(Q)<1$.
Let $P_Q$ be the concentric dilate of $Q$ with $\ell(P_Q)=10\sqrt{n} l(> 10\sqrt{n} l \ell(Q))$.
We write
\begin{equation*}
f_j=f_j\chi_{P_Q}+f_j\chi_{(P_Q)^c}=:\mathtt{f}_j^{(0)}+\mathtt{f}_j^{(1)} \q \text{ for }~j=1,\dots,l.
\end{equation*}
Then 
\begin{align*}
&\Big(\frac{1}{|Q|}\int_Q  \big|T_{\sigma}\big(f_1,\dots,f_l\big)(y)-c_Q\big|^{\frac{r}{l}} \; dy \Big)^{\frac{l}{r}}\\
&\lesssim  \bigg(\frac{1}{|Q|}\int_Q  \Big|T_{\sigma}\big(\mathtt{f}_1^{(0)},\dots,\mathtt{f}_l^{(0)}\big)(y)-d_Q\Big|^{\frac{r}{l}} \; dy \bigg)^{\frac{l}{r}}\\
&\q + \sum_{\substack{\mu_1,\dots,\mu_l\in \{0,1\}\\ (\mu_1,\dots,\mu_l)\not= (0,\dots,0) }}\bigg(\frac{1}{|Q|}\int_Q  \Big|  T_{\sigma}\big(\mathtt{f}_1^{(\mu_1)},\dots,\mathtt{f}_l^{(\mu_l)}\big)(y)-\sum_{k: 2^k \ell(Q)<1}T_{\sigma_k} \big(\mathtt{f}_1^{(\mu_1)},\dots,\mathtt{f}_l^{(\mu_l)}\big)(x)    \Big|^{\frac{r}{l}} \; dy \bigg)^{\frac{l}{r}}
\end{align*}
where
$$d_Q:=c_Q  -\sum_{\substack{\mu_1,\dots,\mu_l\in \{0,1\}\\ (\mu_1,\dots,\mu_l)\not= (0,\dots,0) } }\sum_{k:2^k\ell(Q)<1}T_{\sigma_k}\big(\mathtt{f}_1^{(\mu_1)},\dots,\mathtt{f}_l^{(\mu_l)}\big)(x).    $$
Therefore, in order to establish \eqref{mainestimate}, we need to estimate
\begin{equation*}
 \mathcal{I}_1:=\inf_{d_Q\in \mathbb{C}}\bigg(\frac{1}{|Q|}\int_Q  \Big|T_{\sigma}\big(\mathtt{f}_1^{(0)},\dots,\mathtt{f}_l^{(0)}\big)(y)-d_Q\Big|^{\frac{r}{l}} \; dy \bigg)^{\frac{l}{r}}
\end{equation*}
and for $(\mu_1,\dots,\mu_l)\in \{0,1\}^l\setminus \{(0,\dots,0)\}$
\begin{equation*}
\mathcal{I}_2^{(\mu_1,\dots,\mu_l)}:=\bigg(\frac{1}{|Q|}\int_Q  \Big|  T_{\sigma}\big(\mathtt{f}_1^{(\mu_1)},\dots,\mathtt{f}_l^{(\mu_l)}\big)(y)-\sum_{k: 2^k \ell(Q)<1}T_{\sigma_k} \big(\mathtt{f}_1^{(\mu_1)},\dots,\mathtt{f}_l^{(\mu_l)}\big)(x)    \Big|^{\frac{r}{l}} \; dy \bigg)^{\frac{l}{r}}.
\end{equation*}

Using Proposition \ref{keypropo1}, we have
 \begin{equation*}
 \mathcal{I}_1\lesssim \mathbf{M}_r\big(f_1,\dots,f_l \big)(x).
 \end{equation*}

 For the other terms, we assume $(\mu_1,\dots,\mu_l)\in \{0,1\}^l\setminus \{(0,\dots,0)\}$ and write
   \begin{align*}
&\Big|  T_{\sigma}\big(\mathtt{f}_1^{(\mu_1)},\dots,\mathtt{f}_l^{(\mu_l)}\big)(y)-\sum_{k: 2^k \ell(Q)<1}T_{\sigma_k} \big(\mathtt{f}_1^{(\mu_1)},\dots,\mathtt{f}_l^{(\mu_l)}\big)(x)    \Big|\\
 &\le \sum_{k:2^k\ell(Q)<1}\Big| T_{\sigma_k}\big( \mathtt{f}_1^{(\mu_1)},\dots,\mathtt{f}_l^{(\mu_l)}\big)(y) -T_{\sigma_k}\big( \mathtt{f}_1^{(\mu_1)},\dots,\mathtt{f}_l^{(\mu_l)}\big)(x)    \Big|\\ 
 &\q + \sum_{k:2^k\ell(Q)\ge 1}\Big|  T_{\sigma_k}\big(\mathtt{f}_1^{(\mu_1)},\dots,\mathtt{f}_l^{(\mu_l)}\big)(y)  \Big|\\
 &=:\mathscr{I}_1^{(\mu_1,\dots,\mu_l)}(y)+\mathscr{I}_2^{(\mu_1,\dots,\mu_l)}(y).
 \end{align*}
By applying \eqref{keypropo2_2}, if $y\in Q$, then we have
 $$\Big| T_{\sigma_k}\big( \mathtt{f}_1^{(\mu_1)},\dots,\mathtt{f}_l^{(\mu_l)}\big)(y) -T_{\sigma_k}\big( \mathtt{f}_1^{(\mu_1)},\dots,\mathtt{f}_l^{(\mu_l)}\big)(x)    \Big|\lesssim 2^k\ell(Q)\mathbf{M}_r\big(f_1,\dots,f_l \big)(x) $$
 and thus
 $$\mathscr{I}_1^{(\mu_1,\dots,\mu_l)}(y)\lesssim \mathbf{M}_r\big(f_1,\dots,f_l \big)(x)  \sum_{k:2^k\ell(Q)<1} 2^k\ell(Q)\lesssim \mathbf{M}_r\big(f_1,\dots,f_l \big)(x).  $$
 This yields that
 $$\bigg(\frac{1}{|Q|}\int_Q  \Big(\mathscr{I}_1^{(\mu_1,\dots,\mu_l)}(y) \Big)^{\frac{r}{l}}\; dy\bigg)^{\frac{l}{r}}\lesssim \mathbf{M}_r\big(f_1,\dots,f_l \big)(x).$$
 Furthermore, Proposition \ref{newkeypropo}, together with the property
  $$\big| \mathtt{f}_j^{(\mu_j)}(y) \big|\le \big| f_j(y)\big| \q j=1,\dots,l,$$
 yields 
  \begin{equation*}
\bigg(\frac{1}{|Q|}\int_Q  \Big(\mathscr{I}_2^{(\mu_1,\dots,\mu_l)}(y) \Big)^{\frac{r}{l}}\; dy\bigg)^{\frac{l}{r}}\lesssim \mathbf{M}_r\big(f_1,\dots,f_l \big)(x),
  \end{equation*}
which proves
 \begin{equation*}
 \mathcal{I}_2^{(\mu_1,\dots,\mu_l)}\lesssim \mathbf{M}_r\big(f_1,\dots,f_l \big)(x)
 \end{equation*}
for $(\mu_1,\dots,\mu_l)\in \{0,1\}^l\setminus \{(0,\dots,0)\}$.
This completes the proof.



\end{document}